\documentclass[11pt]{amsart}
\usepackage{amsmath,amssymb,amsthm,amscd}
\usepackage{tabularx}
\usepackage{amsrefs,amsmath,amsfonts,amsthm,amssymb,amscd,comment,xypic,verbatim,cite,xcolor,tabu,hyperref}
\usepackage{longtable}

\numberwithin{equation}{section}

\newtheorem{thm}[equation]{Theorem}

\newtheorem{lemma}[equation]{Lemma}
\newtheorem{prop}[equation]{Proposition}

\theoremstyle{remark}
\newtheorem*{remark}{Remark}
\newtheorem{example}[equation]{Example}

\newcommand{\abs}[1]{\left\lvert#1\right\rvert}

\renewcommand{\bar}[1]{#1\llap{$\overline{\phantom{\rm#1}}$}}

\newcommand{\RR}{\ensuremath{\mathbb{R}}}

\newcommand{\Z}{\ensuremath{\mathbb{Z}}}
\newcommand{\QQ}{\mathbb{Q}}

\newcommand{\QB}{\bar{\mathbb{Q}}}

\DeclareMathOperator{\Pic}{Pic}

\title{Quartic Integral Polynomial Pell Equations}

\author{Zachary Scherr}
\address{
  Google,
  6425 Penn Ave
  Pittsburgh, PA 15206 USA
}
\email{zlscherr@gmail.com}
\author{Katherine Thompson}
\address{
  Department of Mathematics,
  United States Naval Academy,
  121 Blake Road,
  Annapolis, MD 21402 USA
}
\email{kthomps@usna.edu}

\date{\today}

\begin{document}


\begin{abstract}
In this paper we classify all monic, quartic, polynomials $d(x)\in\Z[x]$ for which the Pell equation
\[f(x)^2-d(x)g(x)^2=1\]
has a non-trivial solution with $f(x),g(x)\in\Z[x]$.
\end{abstract}

\maketitle


\section{Introduction}\label{Introduction}
Let $d \geq 2$ be a non-square integer. It is well known that the \textit{Pell equation}
\begin{equation}\label{intro:pell1}
x^2-dy^2 = 1
\end{equation}
has infinitely many solutions $(x,y) \in \mathbb Z^2$ beyond the trivial solutions given by $(x,y) = (\pm 1, 0)$. Moreover there is a \textit{fundamental solution} to these equations: a particular $(x_1,y_1) \in \mathbb N^2$ such that all nontrivial solutions to (\ref{intro:pell1}) can be written as $(\pm x_n, \pm y_n)$ with $n \in \mathbb N$ and $$x_n+y_n\sqrt{d} = (x_1+y_1\sqrt{d})^n.$$ This fundamental solution can be found explicitly for any $d$ by considering the periodicity of the continued fraction expansion of $\sqrt{d}$.
\\
In \cite{Euler}, Euler observed that for special values of $d$, we can relate the fundamental solution of (\ref{intro:pell1}) with a more general polynomial expression. For example, if $d=n^2+1$ for some $n \in \mathbb N$, then the fundamental solution to (\ref{intro:pell1}) is $(x_1,y_1)=(2n^2+1,2n)$ and indeed $$(2n+1)^2-(n^2+1)(2n)^2=1.$$ We can interpret Euler's result as a special case of a more general polynomial Pell's equation. Let $R$ be a commutative ring with unity and let $d(x) \in R[x]$ be a non-square polynomial. We say $d(x)$ is \textit{Pellian over $R[x]$} if there exist $f(x),g(x) \in R[x]$ with $g(x) \neq 0$ and 
\begin{equation}\label{intro:pell2}
f(x)^2-d(x)g(x)^2=1.
\end{equation}
In the polynomial setting, for fixed rings $R$, one can attempt to classify all Pellian polynomials. Abel \cite{Abel} studied the case of $R= \mathbb C$, and his results extend to any field $K$ of characteristic $0$. Similar to the integer case, Abel considered the continued fraction expansion of $\sqrt{d(x)}$, proving that $d(x)$ is Pellian if and only if the continued fraction of $\sqrt{d(x)}$ is periodic.\\
\\
Yet unlike the situation over the integers, as we will see in Section~\ref{Background} the continued fraction of $\sqrt{d(x)}$ need not be periodic. Therefore, in working with the polynomial Pell equation there is an immediate separate question of classifying those $d(x)$ for which (\ref{intro:pell2}) could have a nontrivial solution in the first place.\\
\\
Over the last fifty years several authors have investigated Pellian polynomials over $R= \mathbb Z$. Nathanson \cite{Nathanson} considered $d(x) = x^2+k$, proving that only when $k = \pm 1, \pm 2$ would $d(x)$ be Pellian. Other authors in \cite{Ramasamy} and \cite{Mollin} found further families of Pellian $d(x) \in \mathbb Z[x]$.\\
\\
Perhaps some of the better-known and more recent results, however, come from a series of papers by Webb and Yakota. Instead of looking at particular families of $d(x)$, they considered fixing the period length of the continued fraction expansion of $\sqrt{d(x)}$. In \cite{WY1} they classified the monic $d(x) \in \mathbb Z[x]$, Pellian over $\mathbb Z[x]$, for which the continued fraction expansion of $\sqrt{d(x)}$ has period $1$ or $2$. In particular, since for any monic quadratic polynomial $d(x)$ the continued fraction expansion of $\sqrt{d(x)}$ will have period of length at most $2$, this work answers in full the question of monic quadratic Pellian polynomials.\\
\\
Continuing to generalize, in \cite{Yokota1} Yokota proved that if $d(x) \in \mathbb Z[x]$ is monic of degree at least $4$ and if $d(x)$ is Pellian over $\mathbb Q[x]$, then the period of the continued fraction expansion of $\sqrt{d(x)}$ is even. In \cite{Yokota2} he then proved that if $d(x) \in \mathbb Z[x]$ is monic and quartic with a continued fraction expansion of $\sqrt{d(x)}$ periodic of period $4$ then $d(x)$ is not Pellian over $\mathbb Z[x]$. Additionally, in that paper Yokota asked whether there were Pellian monic, quartic $d(x) \in \mathbb Z[x]$ with continued fraction expansion of $\sqrt{d(x)}$ periodic of period greater than $4$. Our main results answer Yokota's question:
\begin{thm}\label{intro:main}
Let $d(x) \in \mathbb Z[x]$ be a monic, quartic, square-free polynomial for which the continued fraction expansion of $\sqrt{d(x)}$ is periodic of period at least $4$. Then $d(x)$ is Pellian over $\mathbb Z[x]$ if and only if, after a possible $\mathbb Z$-linear change of variables $x\mapsto \pm x + c$ $$d(x) = (x^2-2)(x^2+2x-5).$$
\end{thm}
\begin{thm}\label{nsqf:main}
Let $d(x)\in\Z[x]$ be monic, quartic and nonsquare-free.  Then $d(x)$ is Pellian over $\mathbb Z[x]$ if and only if up to a $\mathbb Z$-linear change of variables $x\mapsto \pm x+c$, $d(x)$ is one of
\begin{itemize}
\item $x^2(x^2\pm 1)$
\item $x^2(x^2\pm 2)$
\item $x^2(x^2-2x-1).$	
\end{itemize}
\end{thm}
\begin{remark}
The $d(x)$ in Theorem~\ref{nsqf:main} all have $\sqrt{d(x)}$ with continued fraction expansion of period at most two.  Thus they are technically covered by Webb and Yokota's work in \cite{WY1}.  That said, in our proof we will additionally show
that there are no other such nonsquare-free $d(x)$ where the continued fraction expansion of $\sqrt{d(x)}$ has period greater than two.
\end{remark}
Whereas most of the previous results involved direct examination and computation of continued fraction expansions, we instead study this problem geometrically via the Jacobian of the curve $$y^2=d(x).$$ In the case $d(x)$ is monic, quartic and square-free, this Jacobian is an elliptic curve.\\
\\
We begin by using Mazur's theorem \cite{Mazur} along with Kubert's parameterization \cite{Kubert} to list all elliptic curves over $\mathbb Q$ with a rational torsion point. Combining this with work of Adams and Razar \cite{AR} we then can list all monic, quartic, square-free $d(x) \in \mathbb Q[x]$ which are Pellian over $\mathbb Q[x]$. Then, we use these parameterizations along with diophantine methods to determine when $f(x), g(x), d(x) \in \mathbb Z[x]$. \\
\\
The paper is organized as follows: in Section 2 we review the theory of continued fractions over Laurent series fields, and the connection between the Pell equation and the geometry of the Jacobian of the curve $y^2=d(x)$. In Section 3 we prove Theorem \ref{nsqf:main} (as the above-mentioned elliptic curve techniques do not apply). In Section 4 we develop parameterizations for all monic, quartic square-free $d(x) \in \mathbb Q[x]$ which are Pellian over $\mathbb Q[x]$. Finally, in Section $5$ we show how to use these parameterizations to prove Theorem \ref{intro:main}.

\section{Background}\label{Background}

Let $K$ be a field with $char(K) \neq 2$ and let $d(x)\in K[x]$.  In this section we investigate when $d(x)$ is Pellian over $K[x]$. Readily apparent is that in order to be Pellian, $\deg(d(x))$ must be even, cannot be a perfect square, and the leading coefficient must be a perfect square in $K^\times$. Without loss of generality we assume $d(x)$ is monic. \\
\\
We begin by defining a continued fraction expansion of $\sqrt{d(x)}$.  To do so we work over the Laurent series field $K((x^{-1}))$.  Any element of this field can be represented by a power series
\[h(x) = \sum_{n\le N}h_nx^n\]
with $h_n\in K$.  For such a series $h(x)$ we define the \textbf{principal part} to be
\[\lfloor h(x)\rfloor = \sum_{n\ge 0} h_nx^n.\]
To construct the continued fraction of $h(x)$ we set
\[\alpha_0(x) = \lfloor h(x)\rfloor\]
and for $n \geq 1$ define
\[\alpha_{n}(x) = \frac{1}{\alpha_{n-1}(x)-\lfloor\alpha_{n-1}(x)\rfloor}.\]
Then the \textbf{continued fraction expansion of $h(x)$} is 
\[[a_0(x);a_1(x),a_2(x),\ldots]\]
where $a_i(x) = \lfloor\alpha_i(x)\rfloor.$\\
\\
Many properties of continued fractions over $\RR$ extend to continued fractions over $K((x^{-1}))$.  For example, if $a \in \mathbb Q$ then the continued fraction expansion of $a$ is finite. Similarly,  $h(x)$ is the Laurent series of a rational function if and only if its continued fraction expansion is finite. \\
\\
Furthermore, if
\begin{align*}
\frac{p_n(x)}{q_n(x)} &= [a_0(x);a_1(x),\ldots,a_n(x)]\\
&= a_0(x) + \cfrac{1}{a_1(x) + \cfrac{1}{\ddots \, + \cfrac{1}{a_n(x)}}}
\end{align*}
then $\frac{p_n(x)}{q_n(x)}$ converges to $h(x)$ in $K((x^{-1}))$ as $n\to\infty.$\\
\\
One of the first results towards solving the polynomial Pell equation was the following:
\begin{thm}[Abel]\label{back:Abel}
Let $d(x)\in K[x]$ be a monic, non-square polynomial of even degree.  The following are equivalent:
\begin{enumerate}
\item $d(x)$ is Pellian over $K[x]$.

\item The continued fraction of $\sqrt{d(x)}$ is periodic.
\item The continued fraction of $\sqrt{d(x)}$ is of the form
\[[a_0(x);\overline{a_1(x),a_2(x),\ldots,a_{n-1}(x),a_n(x)}]\]
where $a_n(x) = 2a_0(x)$ and the sequence
\[a_1(x),a_2(x),\ldots,a_{n-1}(x)\]
is palindromic.
\end{enumerate}
\end{thm}
Theorem \ref{back:Abel} is proved very much in the same way as with the analogous statement over $\Z$. Furthermore, up to constant factors the solutions to the polynomial Pell's equation can be found among the convergents $p_n(x),q_n(x)$ of $\sqrt{d(x)}$. For a proof of Theorem \ref{back:Abel}, see \cite{AR}. 

\begin{example}
Returning to Euler's case where $d(x) = x^2+1$, then
\[\sqrt{x^2+1} = x\sqrt{1+x^{-2}}=x+O(x^{-1}).\]
Thus
\begin{align*}
a_0(x) &= x\\
\alpha_1(x) &= \frac{1}{x\sqrt{1+x^{-2}}-x}\\
&=x\sqrt{1+x^{-2}}+x 	
\end{align*}
yielding the continued fraction
\[\sqrt{x^2+1} = [x;\overline{2x}].\]

Then \[\frac{p_1(x)}{q_1(x)} = \frac{x}{1}\] and we have
\[x^2-(x^2+1)(1)^2 = -1.\]
If $\sqrt{-1}\in K$ then we divide both sides by $\sqrt{-1}$ to get the smallest degree non-trivial solution to the Pell equation.  Otherwise,
\[\frac{p_2(x)}{q_2(x)} = \frac{2x^2+1}{2x}\]
and
\[(2x^2+1)^2-(x^2+1)(2x)^2=1,\] which was precisely Euler's result over $\mathbb Z[x]$.
\end{example}
In general, the smallest $k\ge 1$ for which
\[p_k(x)^2-d(x)q_k(x)^2=c\in K\]
will yield the the smallest degree non-trivial solution to the Pell equation if $\sqrt{c}\in K$.  If not, then
\[p_{2k}(x)+q_{2k}(x)\sqrt{d(x)} = (p_k(x)+q_k(x)\sqrt{d(x)})^2\]
and so
\[p_{2k}(x)^2-d(x)q_{2k}(x)^2 = c^2\]
will yield the smallest degree non-trivial solution to the Pell equation.\\
\\
A major difference between the integer Pell equation and the polynomial Pell equation is that the continued fraction expansion of $\sqrt{d(x)}$ need not be periodic if $K$ is infinite. The proof that for $m \in \mathbb N$, $\sqrt{m}$ has a periodic continued fraction expansion relies on showing the height of $\alpha_n$ is bounded and so eventually must repeat. In the case of $d(x) \in K[x]$, the same logic only shows that the degrees of the numerator and denominator are bounded which is not enough to claim the continued fraction expansion of $\sqrt{d(x)}$ is periodic. In fact, we see later as a consequence of Theorem~\ref{back:geom} that the continued fraction of $\sqrt{d(x)}$ need not be periodic.\\ 
\\
In Theorem~\ref{intro:main} we work over $R=\Z$.  This adds the extra difficulty that we not only require that the continued fraction expansion of $\sqrt{d(x)}$ to be periodic, but we require that some solution to (\ref{intro:pell2}) has integer coefficients.  In the case where the continued fraction expansion of $\sqrt{d(x)}$ has period $1$ or $2$, Webb and Yokota \cite{WY1} were able to bypass this difficulty by proving that if some nontrivial solution to (\ref{intro:pell2}) has $f(x),g(x) \in \mathbb Z[x]$, then necessarily the minimal degree solution lies in $\mathbb Z[x]$.  Thus by explicitly computing the continued fraction expansion of $\sqrt{d(x)}$ one could analyze the coefficients of either $\frac{p_1(x)}{q_1(x)}$ or $\frac{p_2(x)}{q_2(x)}$ to determine when they are integers.  \\
\\
Yet the work of Webb and Yokota concerns monic polynomials. Unfortunately, it need not be the case if $d(x)$ is non-monic that the minimal solution to (\ref{intro:pell2}) have integer coefficients.  In \cite{WY2}, they give examples of this phenomenon.\\ 
\\
In order to assess whether or not the continued fraction of $\sqrt{d(x)}$ is periodic, and equivalently that $d(x)$ is Pellian, we use the geometry of the function field $K(x,y)$ where $y^2 = d(x).$ To do so we make the additional assumption that $d(x)$ is \textit{square-free} and has degree at least $4$.  If $\deg(d(x)) = 2m$ then the function field $K(x,y)$ has genus $(m-1)$.  If we let
$\infty_+$ and $\infty_-$ denote the two infinite places of $K(x,y)$ then one can prove the following:
\begin{thm}\label{back:geom}
Let $d(x)\in K[x]$ be monic, square-free of even degree at least $4$.  If $K(x,y)$ is the function field of $y^2=d(x)$ and $\infty_+$ and $\infty_-$ are the two infinite places of this field over $K(x)$, then
\[f(x)^2-d(x)g(x)^2=1\]
has a non-trivial solution $f(x),g(x)\in K[x]$ if and only if $\infty_+ - \infty_-$ is torsion in $\Pic^0(K(x,y)).$
\end{thm}
\begin{proof}
See \cite{AR}.
\end{proof}
For the proof of Theorem \ref{intro:main}, as $\deg(d(x)) = 4$, $\Pic^0(C)$ is an elliptic curve. Mazur's theorem \cite{Mazur} guarantees that if $P$ is a torsion point on an elliptic curve $E/\QQ$ then the order of $P$ is one of

$$\{1,2,3,4,5,6,7,8,9,10,12\}.$$

In their investigation into the quartic case, Adams and Razar\cite{AR} prove the following: 
\begin{thm}[Adams-Razar]\label{param:arthm1}
Let $d(x)$ be a monic, quartic, square-free polynomial in $K[x]$ and suppose that the continued fraction of $\sqrt{d(x)}$ has period $n$.  If the order of the torsion point $\infty_+-\infty_-$ is $m$ then either
\[n = m-1\]
or
\[n = 2(m-1)\]
where the second case can only occur if $m$ is even.
\end{thm}

Combining this result with aforementioned works of Yokota, Webb and Mazur yields: 
\begin{prop}\label{param:prop1}
Let $d(x)\in\Z[x]$ be monic, square-free and quartic and suppose that $d(x)$ is Pellian over $\mathbb Z[x]$. Let $n$ denote the period of the continued fraction of $\sqrt{d(x)}$ and let $m$ denote the torsion order of the point $\infty_+-\infty_-$ on $\Pic^0(C)$ where $C$ is the curve $y^2=d(x)$.  If $n > 2$ then $(n,m)$ must be one of
\[\{(6, 4), (4, 5), (10, 6), (6, 7), (14, 8), (8, 9), (18, 10), (22, 12)\}.\] 
\end{prop}

Furthermore, work of Kubert~\cite{Kubert} parametrizes all elliptic curves with torsion points over $\QQ$.  Using explicit formulas of Adams and Razar~\cite{AR} allows 
us to parametrize all monic, quartic, square-free $d(x)\in\QQ[x]$ which are Pellian over $\QQ[x]$.  
Within this list of parametrizations we can produce all $d(x)\in\mathbb Z[x]$; however,
unlike in the work of Webb and Yokota we no longer know a priori that the existence of a solution with integer coefficients implies the minimal degree solution has integer coefficients as well.  We will, however, prove that the existence of a solution with integer coefficients necessarily imples that the leading coefficient of the minimal degree solution is an integer.  In Section~\ref{Diophantine} we combine this observation with our parametrizations to complete the proof of Theorem~\ref{intro:main}.

\section{Nonsquare-free $d(x)$}\label{Squarefree}

The geometric methods we use in proving Theorem~\ref{intro:main} only apply when $d(x)$ is square-free and the period of the continued fraction expansion of $\sqrt{d(x)}$ is at least three.  In this section we analyze the case of $d(x)\in\Z[x]$ monic, quartic and nonsquare-free and prove Theorem \ref{nsqf:main}.

We begin by noting that any monic, quartic, non-squarefree $d(x)\in\Z[x]$ can be written in the form
\[d(x) = D(x)(x-a)^2\]
where $D(x)\in\Z[x]$ is monic and quadratic, and $a\in\Z$.  Thus any non-trivial solution to the Pell equation
\[f(x)^2-d(x)g(x)^2=1\]
with $f(x),g(x)\in\Z$ can be rewritten as
\[f(x)^2-D(x)((x-a)g(x))^2=1.\]
This observation gives rise to the following:
\begin{lemma}\label{nsqf:lemma1}
The monic, quartic, non-squarefree $d(x)\in\Z[x]$ which are Pellian over $\mathbb Z[x]$ are in correspondence with monic, quadratic, Pellian $D(x) \in \mathbb Z[x]$ for which
\[F(x)^2-D(x)G(x)^2=1\]
where $G(x)$ has a root in $\Z$.
\end{lemma}

We also make use of the following propositions, standard in the theory of the Pell equation (see \cite{schmidt2000continued} for a detailed introduction to polynomial Pell equations):
\begin{prop}\label{nsqf:prop1}
Let $K$ be a field and suppose $d(x)\in K[x]$ is Pellian over $K[x]$. If $f_1(x), g_1(x)$ give such a solution where $\deg(f_1(x))$ is minimal, 
then every solution to (\ref{intro:pell2}) is given by $\pm f_n(x), \pm g_n(x)$ where
\[f_n(x)+g_n(x)\sqrt{d(x)} = (f_1(x)+g_1(x)\sqrt{d(x)})^n\] and $n \in \mathbb Z_{\geq 0}$. 
\end{prop}
\begin{prop}\label{nsqf:prop2}
Let $K, d(x), f_n(x), g_n(x)$ be as in Proposition~\ref{nsqf:prop1}.  Then the $g_n(x)$ satisfy
the linear recurrence
\[g_{n+2}(x) = 2f_1(x)g_{n+1}(x) - g_n(x).\]

Furthermore, if $d(x)$ is monic and $c_n\in K$ is the leading coefficient of $g_n(x)$, then
\[c_n = 2^{n-1}c_1^n.\]
\end{prop}

\begin{proof}[Proof of Theorem~\ref{nsqf:main}]

According to Lemma~\ref{nsqf:lemma1}
we let $D(x)=x^2+rx+s\in\Z[x]$.  Making the substitution $x\mapsto x-\lfloor\frac{r}{2}\rfloor$ allows
us to assume without loss of generality that either
\[D(x) = x^2+s\]
or
\[D(x) = x^2+x+s\]
for some $s\in\Z$.

If $D(x)=x^2+s$, then the continued fraction of $\sqrt{D(x)}$ in $\QQ((x^{-1}))$ is
\[\sqrt{D(x)} = \left[ x;\overline{\frac{2x}{s},2x}\right].\]
Thus the period of the continued fraction of $\sqrt{D(x)}$ is one if $s=1$ and two otherwise.  Using the continued fraction for
$\sqrt{D(x)}$ we find that the smallest degree non-trivial solution to (\ref{intro:pell2}) over $\QB$ is
\[
(F_1(x),G_1(x)) = \left(\frac{1}{\sqrt{-s}}x,\frac{1}{\sqrt{-s}}\right).
\]
Suppose that some solution $F_n(x),G_n(x)\in\Z[x]$.  Then by Proposition~\ref{nsqf:prop1} we know that $F_{2n}(x),G_{2n}(x)\in\Z[x]$
as well.  The leading coefficient of $G_{2n}(x)$ is
\[\frac{2^{2n-1}}{(-s)^n}\in\Z\]
which forces $s\in\{\pm 1,\pm 2\}.$

If $s = -1$ then the smallest degree solution to the Pell equation over $\QQ$ is $(F_1(x),G_1(x)) = (x,1)$.  If, on the other hand,
$s \in \{1, \pm 2\}$ then over $\QQ$ the smallest degree solution to the Pell equation is
\[(F_1(x),G_1(x)) = \left( \frac{2}{s}x^2+1,\frac{2}{s}x \right).\]
Thus for any $s\in\{\pm 1,\pm 2\}$, the smallest solution to (\ref{intro:pell2}) for $D(x) =x^2+s$ over $\mathbb Q[x]$ has integer coefficients, and hence
every solution has integer coefficients.  Therefore, we need only find all $a\in\Z$ and $n\ge 1$
so that $G_n(a) = 0$ in accordance with Lemma~\ref{nsqf:lemma1}.

To this end we can use Proposition~\ref{nsqf:prop2} to deduce that the sequence $G_n(a)$ satisfies the recurrence relation
\[G_{n+2}(a) = 2F_1(a)G_{n+1}(a) - G_n(a).\]
We now proceed by cases on $s$. 

{\bf Case $s=-1$:} 

Here $F_1(x) = x$ and $G_1(x) = 1$ so $G_0(a)=0, G_1(a) = 1$ and

$$G_{n+2}(a) = 2aG_{n+1}(a)-G_n(a).$$

Dividing by $G_{n+1}(a)$ and taking absolute values yields

$$\abs{\frac{G_{n+2}(a)}{G_{n+1}(a)}} = \abs{2a - \frac{G_n(a)}{G_{n+1}(a)}}\ge 2\abs{a}-\abs{\frac{G_n(a)}{G_{n+1}(a)}}. $$

If $a\ne 0$ then an induction argument shows that $\abs{G_{n+1}(a)} > \abs{G_n(a)}$ for $n\ge 0$, and thus $G_n(a) \neq 0$.  
If $a = 0$ then $G_2(a) = 0$, giving rise through Lemma~\ref{nsqf:lemma1} to $d(x) = x^2(x^2-1)$.

{\bf Case $s = 1, \pm 2$:} 

Here $F_1(x) = \frac{2}{s}x^2+1$ and $G_1(x) = \frac{2}{s}x$.  Thus $G_0(a) = 0, G_1(a) = \frac{2a}{s}$ and

$$G_{n+2}(a) = 2\left(\frac{2}{s}a^2+1\right)G_{n+1}(a)-G_n(a).$$

Again we divide by $G_{n+1}(a)$ and take absolute values:

$$\abs{\frac{G_{n+2}(a)}{G_{n+1}(a)}} \ge 2\abs{\frac{2}{s}a^2+1}-\abs{\frac{G_n(a)}{G_{n+1}(a)}}.	$$

If $a = 0$, $G_1(a) = 0$ and $d(x) = x^2(x^2+1)$ or $d(x)=x^2(x^2\pm2)$.  So we suppose that $a\ne 0$ and so in particular $G_1(a)\ne 0$.  If $\frac{2}{s}a^2+1 \ne 0$ then since it is an integer we have $\abs{\frac{2}{s}a^2+1} \ge 1$ and again a simple induction argument shows $\abs{G_{n+1}(a)} > \abs{G_n(a)}$ and hence $G_n(a)\ne 0$ for any $n\ge 1$.

Thus we only need to consider $\frac{2}{s}a^2+1=0$ which occurs if $s=-2$ and $a=\pm 1$.  In this case $G_2(\pm 1) = 0$ and so
\[d(x) = (x\pm 1)^2(x^2-2).\]
Making the substitution $x\mapsto \pm x\pm 1$ transforms these two possibilities into
\[d(x) = x^2(x^2-2x-1).\]

Thus we have so far shown that the case $D(x) = x^2+s$ yields all the $d(x)$ in Theorem~\ref{nsqf:main}.  We now show that the case $D(x)=x^2+x+s$ where $s\in\Z$ yields no additional $d(x)$.

Working over $\QQ((x^{-1}))$ we have
\[\sqrt{D(x)} = [x+\frac{1}{2};\overline{\frac{2x}{s-1/4}+\frac{1}{s-1/4}, 2x+1}].\]
Thus over $\QB$, the smallest degree non-trivial solution to the Pell equation $F(x)^2-D(x)G(x)^2=1$ is
\[(F_1(x),G_1(x)) = \frac{1}{\sqrt{1/4-s}}(x+1/2,1).\]
If we have a solution $F_n(x),G_n(x)\in\Z[x]$ then we have $F_{2n}(x),G_{2n}(x)\in\Z[x]$. The leading coefficient of
$G_{2n}(x)$ is 
\[\frac{2^{2n-1}}{(1/4-s)^n} = \frac{2^{4n-1}}{(1-4s)^n}
\]
and so $G_{2n}(x)\in\Z[x]$ implies $s = 0$.  Thus $D(x)=x^2+x$ and the minimal non-trivial solution
to the Pell equation is
\[(F_1(x),G_1(x)) = (2x+1,2).\]
Similar techniques as in the previous case will show that $G_n(a)\ne 0$ for any $a\in\Z$ and $n\ge 1$, completing the proof.  
\end{proof}

\section{Parametrizations}\label{Parametrizations}
Over $\QQ$ Kubert\cite{Kubert} gives (up to isomorphism) parameterizations of all elliptic curves $E$ with a torsion point $P\in E(\QQ)$ of fixed order.  We begin with an elliptic curve $E$ in Tate normal form
\[y^2 + (1-c)xy-by = x^3-bx^2\]
where $P\mapsto(0,0)$.  The discriminant $\Delta$ of this elliptic curve is
\[\Delta(b,c) = b^3(16b^2 -8bc^2 -20bc+b+c^4 -3c^3 +3c^2 -c).\]
Provided $\Delta(b,c)\ne 0$, the following table gives parameterizations for $b$ and $c$ in terms of a parameter $t\in\QQ$ so that $(0,0)$ has order $m$:

\renewcommand{\arraystretch}{1.3}
\begin{table}[ht]
\centering
\caption{Parameterizations}
\begin{tabular}{|c|c|c|}
\hline
$m$&$b$&$c$\\\hline
4&$t$&0\\\hline
5&$t$&$t$\\\hline
6&$t(t+1)$&$t$\\\hline
7&$t^2(t-1)$&$t(t-1)$\\\hline
8&$(2t-1)(t-1)$&$\frac{(2t-1)(t-1)}{t}$\\\hline
9&$t^2(t-1)(t^2-t+1)$&$t^2(t-1)$\\\hline
10&$\frac{t^3(t-1)(2t-1)}{(t^2-3t+1)^2}$&$\frac{-t(t-1)(2t-1)}{t^2-3t+1}$\\\hline
12&$\frac{t(2t-1)(3t^2-3t+1)(2t^2-2t+1)}{(t-1)^4}$&$\frac{-t(2t-1)(3t^2-3t+1)}{(t-1)^2}$\\\hline
\end{tabular}
\end{table}

We use this to parameterize all $d(x)\in\QQ[x]$ for which $d(x)$ is monic, quartic and square-free and Pellian over $\mathbb Q[x]$. First observe that after applying a suitable transformation $x\mapsto x+a$ with $a\in\QQ$ we may assume that $d(x)$ is of the form
\[d(x) = x^4+r_2x^2+r_1x+r_0.\]
To determine $r_2,r_1$ and $r_0$ we begin with the following proposition of Adams and Razar\cite{AR} which gives a concrete model in short Weierstrass for the Jacobian of $y^2=d(x)$.
\begin{prop}[Adams, Razar]\label{param:prop2}
Let $K$ be a field with $char(K) \neq 2,3$. Let
\[d(x) = x^4 - 6ax^2 - 8bx+c\]
be monic and square-free with $a,b\in K$. Then the Jacobian of the curve $y^2=d(x)$ is isomorphic to
\[y^2=x^3+Ax+B\]
where
\begin{align*}
A&=-\tfrac{1}{4}(c+3a^2)\\
B&=b^2-a^3-Aa.
\end{align*}
Moreover, under this isomorphism the point $\infty_+-\infty_-$ maps to $(a,b)$.
\end{prop}

We now begin with Kubert's list of elliptic curves in Tate normal form having a torsion point of order $m$.  We then transform these curve into short Weierstrass form, tracking the image of $(0,0)$.  Noting that Kubert's parametrization is only given up to isomorphism and recalling that over $\QQ$, two elliptic curves $y^2=x^3+rx+s$ and $y^2=x^3+r'x+s'$ are isomorphic iff there is a $u\in\QQ^\times$ so that $r' = \frac{r}{u^4}$ and $s' = \frac{s}{u^6}$, and using Proposition~\ref{param:prop2} allows us to parametrize all of our quartics.  The result we obtain is
\begin{thm}\label{param:thm1}
Let $d(x)=x^4+r_2x^2+r_1+r_0\in\QQ[x]$ be square-free.  Then $d(x)$ is Pellian over $\mathbb Q[x]$ and the torsion order of $\infty_+-\infty_-$ on the Jacobian of $y^2=d(x)$ is $m\ge 4$ if there are $a,b\in\QQ$ so that one of the following is true.
\renewcommand{\arraystretch}{1.3}
\begin{longtable}{|c|| l |}
\hline
$m$ &$ r_i$ \\
\hline
$4$ & \begin{tabular}{l}
$r_2=(8a-2)/b^2$\\
$r_1 = 32a/b^3$ \\
$r_0=(a^4-12a^3+6a^2+20a+1)/b^4$
\end{tabular}
\\
\hline
$5$ & \begin{tabular}{l}
$r_2 = (-2a^2+12a-2)/b^2$\\
$r_1 = 32a/b^3$\\
$r_0=(a^4 - 12a^3 + 6a^2 + 20a + 1)/b^4$
\end{tabular}
\\
\hline
$6$ & \begin{tabular}{l}
$r_2=(6a^2 + 12a - 2)/b^2$\\
$r_1=(32a^2 + 32a)/b^3$\\
$r_0=(9a^4 + 4a^3 + 30a^2 + 20a + 1)/b^4$
\end{tabular}
\\
\hline
$7$ & \begin{tabular}{l}
$r_2=(-2a^4 + 12a^3 - 6a^2 - 4a - 2)/b^2$\\
$r_1=(32a^3 - 32a^2)/b^3$\\
$r_0=(a^8 - 12a^7 + 42a^6 - 64a^5 + 51a^4 - 22a^2 + 4a + 1)/b^4$
\end{tabular}
\\
\hline
$8$ & \begin{tabular}{l}
$r_2=(8a^4 + 8a^3 - 32a^2 + 16a - 2)/(a^2b^2)$\\
$r_1=(64a^2 - 96a + 32)/b^3$
$r_0=(16a^8 - 96a^7 + 336a^6 - 576a^5 + 536a^4 - 296a^3 + 96a^2 - 16a + 1)/(a^4b^4)$
\end{tabular}
\\
\hline
$9$ & \begin{tabular}{l}
$r_2=(-2a^6 + 12a^5 - 18a^4 + 20a^3 - 12a^2 - 2)/b^2$\\
$r_1=(32a^5 - 64a^4 + 64a^3 - 32a^2)/b^3$\\
$r_0=(a^{12} - 12a^{11} + 54a^{10} - 128a^9 + 181a^8 - 156a^7 + 82a^6 - 4a^5 - 42a^4 + 44a^3 - 20a^2 + 1)/b^4$
\end{tabular}
\\
\hline
$10$ & \begin{tabular}{l}
$r_2= (-8a^6 + 32a^5 - 16a^4 - 16a^3 + 8a - 2)/(a^4b^2 - 6a^3b^2 + 11a^2b^2 - 6ab^2 + b^2)$\\
$r_1= (64a^5 - 96a^4 + 32a^3)/(a^4b^3 - 6a^3b^3 + 11a^2b^3 - 6ab^3 + b^3)$\\
$r_0=$  \tiny{$(16a^{12} - 128a^{11} + 448a^{10} - 896a^9 + 1024a^8 - 416a^7 - 408a^6 + 
    608a^5 - 304a^4 + 48a^3 + 16a^2$}\\
  \hspace{.5in}\tiny{ $8a + 1)/(a^8b^4 - 12a^7b^4 + 
    58a^6b^4 - 144a^5b^4 + 195a^4b^4 - 144a^3b^4 + 58a^2b^4 - 12ab^4
    + b^4$)
}
\end{tabular}
\\
\hline
$12$ & \begin{tabular}{l}
$r_2 = $\tiny{$(24a^8 - 240a^7 + 672a^6 - 936a^5 + 744a^4 - 336a^3 + 72a^2 - 2)/(a^6b^2
    - 6a^5b^2 + 15a^4b^2$}\\
     \hspace{.5in}\tiny{$- 20a^3b^2 + 15a^2b^2- 6ab^2 + b^2)$}\\
$r_1 = (384a^6 - 960a^5 + 1088a^4 - 672a^3 + 224a^2 - 32a)/(a^4b^3 - 4a^3b^3 +
    6a^2b^3 - 4ab^3 + b^3)
$\\
$r_0 =$ \tiny{$ (144a^{16} - 576a^{15} + 2112a^{14} - 9696a^{13} + 34016a^{12} - 82176a^{11} + 
    141936a^{10} - 181984a^9$}\\
     \hspace{.5in}\tiny{$+ 177240a^8 - 132528a^7 + 76096a^6 - 33208a^5 +
    10760a^4 - 2480a^3 + 376a^2 - 32a+1)/$}\\
     \hspace{.5in}\tiny{$(a^{12}b^4 - 12a^{11}b^4 + 
    66a^{10}b^4 - 220a^9b^4 + 495a^8b^4 - 792a^7b^4 + 924a^6b^4 - 
    792a^5b^4+ 495a^4b^4 - 220a^3b^4 + 66a^2b^4 - 12ab^4 + b^4)$}
\end{tabular}\\
\hline

\end{longtable}

\end{thm}
\section{Main Proof}\label{Diophantine}

Having parameterized all $d(x) = x^4 + r_2x^2+r_1x+r_0\in\QQ[x]$ that are Pellian over $\mathbb Q[x]$, we now add the constraint of being Pellian over $\mathbb Z[x]$. There are two main hurdles to overcome.  First, the $d(x)$ in our parameterizations lack a cubic term and over $\Z$ we cannot always make such a substitution. 
If $f(x)^2-d(x)g(x)^2 = 1$ with $d(x),f(x),g(x)\in\Z[x]$ then we can find $c\in\QQ$ for which $d(x+c)$ has no cubic term. 
While we \textit{cannot} guarantee that $f(x+c), g(x+c)\in\Z[x]$, we can guarantee that the leading coefficients of $f(x+c)$ and $g(x+c)$ are integers.  As it turns out this will be enough to prove Theorem \ref{intro:main}.

The second hurdle is that with our parameterizations from Theorem~\ref{param:thm1} there are infinitely many $f(x),g(x)\in\QQ[x]$ for which $f(x)^2-d(x)g(x)^2=1$.  These solutions are generated by a minimal degree solution, but it is possible that that minimal degree solution is not in $\Z[x]$.  To this end, we will prove that if $f(x)^2-d(x)g(x)^2=1$ has a solution $f(x),g(x)\in\Z[x]$ then the minimal degree solution has an integer leading coefficient.  

We begin with the following lemmas:
\begin{lemma}\label{main:lemma1}
Let $d(x) = x^4 + a_3x^3 + a_2x^2+a_1x+a_0\in\Z[x]$ and define
\[d_c(x) := d\left(x-\frac{a_3}{4}\right) = x^4+b_2x^2+b_1x+b_0. \]
Then $8b_2, 8b_1, 256b_0\in\Z$.
\end{lemma}
\begin{proof}
Left to reader. 
\end{proof}
\begin{lemma}\label{main:lemma2}
	Let $d(x)\in\Z[x]$ be monic, quartic and square-free.  Suppose $d(x)$ is Pellian over $\mathbb Q[x]$ with some solution $f(x),g(x) \in \mathbb Q[x]$ having integral leading coefficients. Then all solutions have integral leading coefficients.
\end{lemma}
\begin{proof}
Let $f_1(x), g_1(x)$ be a nontrivial solution to \eqref{intro:pell2} where the degree of $f_1(x)$ is minimal.  Then we know that all nontrivial solutions 
are given by $\pm f_n(x), \pm g_n(x)$ where
\begin{equation}\label{main:eqn2}
f_n(x) + g_n(x)\sqrt{d(x)} = (f_1(x) + g_1(x)\sqrt{d(x)})^n
\end{equation}
with $n \ge 1$.  Let $a_n$ denote the leading coefficient of $f_n(x)$.  Since $d(x)$ is monic we know that $a_n$ is also the leading coefficient
of $g_n(x)$.  If we consider \eqref{main:eqn2} in the power series ring $\QQ((x^{-1}))$ then 
\[\sqrt{d(x)} = x^2 + (\text{lower order terms})\]
and so the leading coefficient of the left hand side of \eqref{main:eqn2} is $2a_n$.  The binomial theorem says the leading coefficient of the right hand side is 
\[\sum_{k=0}^n \binom{n}{k}a_1^ka_1^{n-k} = 2^na_1^n.\]
This gives
\begin{equation}\label{main:eqn3}
a_n = 2^{n-1}a_1^n.
\end{equation}
By assumption $a_n = (2a_1)^{n-1}a_1\in\Z$ for some $n\ge 1$.  If for some prime $p$ we have $\abs{a_1}_p > 1$ then
\[\abs{a_n}_p = (\abs{2}_p\abs{a_1}_p)^{n-1}\abs{a_1}_p \ge \abs{a_1}_p > 1\]
contradicting the fact that $a_n\in\Z$.  Thus if any $a_n\in\Z$ we have shown that $a_1\in\Z$.
\end{proof}

We now have:
\begin{prop}\label{main:prop1}
	Let $d(x)\in\Z[x]$ be a monic, quartic, square-free polynomial which is Pellian over $\mathbb Z[x]$. Let $d_c(x)$ be as in Lemma \ref{main:lemma1}. If $f_1(x)^2-d_c(x)g_1(x)^2=1$ is the minimal degree solution to the Pell equation then the leading coefficient of $f_1(x)$ is an integer.
\end{prop}
\begin{proof}
Let $f(x)^2-d(x)g(x)^2=1$ be a non-trivial solution to the Pell equation with $f(x),g(x)\in\mathbb Z[x]$.
Since $d_c(x) = d(x+r)$ for some $r\in\QQ$, we know that
\[f(x+r)^2 - d_c(x)g(x+r)^2 = 1.\]
But the leading coefficient of $f(x+r)$ is the same as the leading coefficient of $f(x)$ which is an integer.  Thus by Lemma~\ref{main:lemma2}, the leading coefficient of $f_1(x)$
must also be an integer.
\end{proof}

Proposition~\ref{main:prop1} suggests the strategy we may follow to prove Theorem~\ref{intro:main}.  We start with our parameterizations from Theorem~\ref{param:thm1}.
For each $d(x) = x^4+r_2x^2+r_1x+r_0$ we compute a non-trivial $f(x),g(x)$ solution to \eqref{intro:pell2}.  By Lemma~\ref{main:lemma1} and Proposition~\ref{main:prop1}, in order
for there to exist $c\in\QQ$ with $d(x+c),f(x+c),g(x+c)\in\Z[x]$ we necessarily need that $8r_2, 8r_1, 256r_0$ and the leading coefficient of $f(x)$ to all be integers.  It turns out these integrality conditions will be enough to restrict the $a,b$ in the parametrizations in Theorem~\ref{param:thm1} to a finite set, from which we can deduce Theorem~\ref{intro:main}.  

As an illustration, consider the case when the torsion order is $4$ and
\[d(x) = x^4+\frac{8a-2}{b^2}x^2 + \frac{32a}{b^3}x + \frac{16a^2+24a+1}{b^4}.\]
Using Sage code \cite{SageCode} we obtain the following continued fraction expansion of $\sqrt{d(x)}$:

\[\left[ x^2 + \frac{4a - 1}{b^2};
  \overline{\frac{b^3}{16a}x+\frac{-b^2}{16a},
  \frac{4}{b}x+\frac{-4}{b^2},
 \frac{b^4}{32a}x^2+\frac{4ab^2 - b^2}{32a},
 \frac{4}{b}x+\frac{-4}{b^2},
 \frac{b^3}{16a}x + \frac{-b^2}{16a},
 2x^2+\frac{8a - 2}{b^2}}\right].\]
By Theorem~\ref{param:arthm1} the period of the continued fraction of $\sqrt{d(x)}$ is either $3$ or $6$, and we see the period is $3$ precisely when $\frac{b^4}{32a} = 2$.

We the compute a solution to the Pell equation
\[f(x)^2-d(x)g(x)^2 = 1,\]
where
\begin{multline*}
f(x)=\frac{256 a^{4} + 224 a^{2} + 32 a + 1}{512 a^{3}} + \left(\frac{64 a^{3} b - 48 a^{2} b - 20 a b - b}{128 a^{3}}\right) x + \\
+ \left(\frac{64 a^{3} b^{2} - 48 a^{2} b^{2} + 4 a b^{2} + b^{2}}{128 a^{3}}\right) x^{2} +
\left(\frac{16 a^{2} b^{3} + 24 a b^{3} + b^{3}}{128 a^{3}}\right) x^{3} + \left(\frac{48 a^{2} b^{4} - 32 a b^{4} - 5 b^{4}}{256 a^{3}}\right) x^{4} + \\
+ \left(\frac{-4 a b^{5} + b^{5}}{128 a^{3}}\right) x^{5} + \left(\frac{4 a b^{6} + b^{6}}{128 a^{3}}\right) x^{6} + \left(\frac{-b^{7}}{128 a^{3}}\right) x^{7} + \frac{b^{8}}{512 a^{3}} x^{8}
\end{multline*}
and
\begin{multline*}
g(x) = \frac{64 a^{3} b^{2} - 48 a^{2} b^{2} - 20 a b^{2} - b^{2}}{512 a^{3}} + \left(\frac{12 a b^{3} + b^{3}}{128 a^{3}}\right) x + \left(\frac{48 a^{2} b^{4} - 24 a b^{4} - 5 b^{4}}{512 a^{3}}\right) x^{2} + \left(\frac{-b^{5}}{32 a^{2}}\right) x^{3} + \\ +\left(\frac{12 a b^{6} + 5 b^{6}}{512 a^{3}}\right) x^{4} + \left(\frac{-b^{7}}{128 a^{3}}\right) x^{5} + \frac{b^{8}}{512 a^{3}} x^{6}.
\end{multline*}

By Lemma~\ref{main:lemma1} and Proposition~\ref{main:prop1}, in order for this solution to correspond to one which is Pellian over $\mathbb Z[x]$ we need 
\[\frac{b^{8}}{512 a^{3}}, 8\left(\frac{8a-2}{b^2}\right), 8\left(\frac{32a}{b^3}\right), 256\left(\frac{16a^2+24a+1}{b^4}\right)\in\Z.\]
In particular, 

\begin{equation}\label{main:eqn4}
\frac{b^{8}}{512 a^{3}}\cdot\left(256\left(\frac{16a^2+24a+1}{b^4}\right)\right)^2 = \frac{32768 a^{4} + 98304 a^{3} + 77824 a^{2} + 6144 a + 128}{a^{3}}\in\Z
\end{equation}

The following lemma shows us how to compute all $a\in\QQ$ for which \eqref{main:eqn4} holds:

\begin{lemma}\label{main:lemma3}
Let $p(x),q(x)\in\Z[x]$ with $\deg p(x) > \deg q(x)$ and $q(0) = 0$. Suppose that $a\in\QQ$ is such that
\[\frac{p(a)}{q(a)}\in\Z.\]
Writing $a=\frac{r}{s}$ in lowest terms and $p(x) = \sum_{i=0}^n p_ix^i$ with $p_n\ne 0$, then $s \vert p_n$ and $r \vert p_0$.
\end{lemma}
\begin{proof}
This is a simple application of the rational root theorem, left to the reader.  
\end{proof}

Applying Lemma~\ref{main:lemma3} to \eqref{main:eqn4} we immediately deduce that $\abs{a}_p = 1$ for $p\ne 2$ and
\[\frac{1}{128}\le\abs{a}_2\le 32768.\]

We use these $p$-adic inequalities to bound the $p$-adic absolute value of $b$.  Since $\frac{b^{8}}{512 a^{3}}\in\Z$ we know 
\[\abs{b}_p^8\le \abs{512}_p\abs{a}_p^3\]
for all primes $p$.  Similarly, since $8\left(\frac{32a}{b^3}\right)\in\Z$ we see
\[\abs{b}_p^3\ge\abs{256}_p\abs{a}_p.\]
These two inequalities yield $\abs{b}_p = 1$ for $p\ne 2$ and \[\left(\frac{1}{32768}\right)^{1/3} \le \abs{b}_2\le (68719476736)^{1/8}.\]
This yields at most $84$ possible pairs $(a,b)\in\QQ^2$ for which 
\[\frac{b^{8}}{512 a^{3}}, 8\left(\frac{8a-2}{b^2}\right), 8\left(\frac{32a}{b^3}\right), 256\left(\frac{16a^2+24a+1}{b^4}\right)\in\Z.\]

Continuing with our example, suppose we look at the pair $(a,b) = (2,4)$.  We have
\[d(x) = x^{4} + \frac{7}{8} x^{2} + x + \frac{113}{256}.\]
By Lemma~\ref{main:lemma1}, we now solve for $c\in\Z$ so that $d\left(x+\frac{c}{4}\right)\in\Z[x]$: 
\[d\left(x+\frac{c}{4}\right) = x^{4} + c x^{3} + \left(\frac{3}{8} c^{2} + \frac{7}{8}\right) x^{2} + \left(\frac{1}{16} c^{3} + \frac{7}{16} c + 1\right) x + \frac{1}{256} c^{4} + \frac{7}{128} c^{2} + \frac{1}{4} c + \frac{113}{256}.\]
However, if $d(x+\frac{c}{4})\in\mathbb Z[x]$ then for any $k\in\mathbb Z$ we would have $d(x+\frac{c+4k}{4})\in\mathbb Z[x]$ as well.  Thus we only need to check $c=0,1,2,3$ and since none of these values of $c$ yield $d(x+\frac{c}{4})\in\mathbb Z[x]$ we have eliminated the pair $(a,b)=(2,4)$.

Repeating this process for the remaining pairs 
reduces our list of $(a,b)$ to 
\[(a,b) = \left(-\frac{1}{16},\pm1\right),\left(-\frac{1}{64},\pm\frac{1}{2}\right).\]

We can eliminate $(a,b) = \left(-\frac{1}{16},\pm 1\right)$ since the resulting $d(x)$ is not square-free.  Finally, for $(a,b) = \left(-\frac{1}{64},\pm\frac{1}{2}\right)$
The corresponding polynomials $d\left(x+\frac{c}{4}\right)$ are
\begin{eqnarray*}
d_1 \left( x+ \frac{c}{4} \right) & = & x^{4} + c x^{3} + \left(\frac{3}{8} c^{2} - \frac{17}{2}\right) x^{2} + \left(\frac{1}{16} c^{3} - \frac{17}{4} c + 4\right) x + \frac{1}{256} c^{4} - \frac{17}{32} c^{2} + c + \frac{161}{16}\\
d_2 \left( x+ \frac{c}{4} \right) & = & x^{4} + c x^{3} + \left(\frac{3}{8} c^{2} - \frac{17}{2}\right) x^{2} + \left(\frac{1}{16} c^{3} - \frac{17}{4} c - 4\right) x + \frac{1}{256} c^{4} - \frac{17}{32} c^{2} - c + \frac{161}{16}.
\end{eqnarray*}
With the substitutions $x\mapsto -x$ and $c\mapsto -c$ we reduce this list further to just:
$$
d\left( x+ \frac{c}{4} \right)  =  x^{4} + c x^{3} + \left(\frac{3}{8} c^{2} - \frac{17}{2}\right) x^{2} + \left(\frac{1}{16} c^{3} - \frac{17}{4} c + 4\right) x + \frac{1}{256} c^{4} - \frac{17}{32} c^{2} + c + \frac{161}{16}
$$
The only $c\in\{0,1,2,3\}$ for which $d(x+\frac{c}{4})\in\mathbb Z[x]$ is $c=2$, and making this substitution yields
$$
d(x)  =  x^{4} + 2 x^{3} - 7 x^{2} - 4 x + 10.
$$

So far we've proved that up to a change of variables $x\mapsto \pm x + c$ with $c\in\mathbb Z$, the polynomial $d(x) = x^{4} + 2 x^{3} - 7 x^{2} - 4 x + 10$ is the unique monic, quartic,
square-free $d(x)\in\mathbb Z[x]$ for which the point $\infty_+ - \infty_-$ on the Jacobian of the curve $y^2=d(x)$ has torsion order $4$ and for which the minimal degree non-trivial solution to
$$f(x)^2-d(x)g(x)^2=1$$
with $f(x),g(x)\in\QQ[x]$ have integral leading coefficients.

Using our continued fraction algorithm we get
\begin{align*}
f(x)&=2x^{8} + 24x^{7} + 100x^{6} + 120x^{5} - 266x^{4} - 792x^{3} - 244x^{2} + 912x + 721\\
g(x)&=2x^{6} + 22x^{5} + 86x^{4} + 118x^{3} - 74x^{2} - 334x - 228.
\end{align*}
and so indeed 
$$d(x) = x^{4} + 2 x^{3} - 7 x^{2} - 4 x + 10 = (x^2 - 2)(x^2 + 2x - 5)$$ is Pellian over $\mathbb Z[x]$.
To complete the proof of Theorem~\ref{intro:main} we must repeat the above analysis for each of the parameterizations from Section~\ref{Parametrizations} to show that no additional $d(x)$ exist.  In fact, what one finds is that for higher torsion order there are no monic, quartic, square-free $d(x)\in\mathbb Z[x]$ for which the minimal degree solution to the Pell equation, $f_1(x),g_1(x)\in\QQ[x]$, has an integral leading coefficient.  Lemma~\ref{main:lemma1}, Proposition~\ref{main:prop1} and Lemma~\ref{main:lemma3} along with the parametrizations in Section~\ref{Parametrizations} turn this problem into a finite computation, and we have provided Sage code \cite{SageCode} detailing this computation for all torsion orders.


\begin{thebibliography}{fszw90}
\bibitem{Abel} Abel, Niels H., \textit{\"Uber die integration der differential-formel, wenn $r$ und $\rho$ ganze functioninen sind}, Journal f\"ur die Reine und Angewandte Mathematik (1826), No. 1, 185-221.
\bibitem{AR} Adams, William W., Razar, Michael J., \textit{Multiples of points on elliptic curves and continued fractions}, Proceedings of the London Mathematica Society \textbf{3} (1980), No. 3, 481-498.
\bibitem{Euler} Euler, Leonhard, \textit{De usu novi algorithmi in problemate pelliano solvendo}, Novi comentarii academiae scientiarum Petropolitanae (1767), 28-66.
\bibitem{Kubert} Kubert, Daniel Sion, \textit{Universal bounds on the torsion of elliptic curves}, Proceedings of the London Mathematical Society \textbf{3} (1976), No. 2, 193-237.
\bibitem{Mazur} Mazur, Barry, \textit{Rational points on modular curves}, Modular Functions of One Variable V (Proc. Internat. Conference, Univ. Bonn, Bonn, 1977), Lecture Notes in Math. \textbf{601}, Springer-Verlag, Berlin (1977), 107-148.
\bibitem{Mollin} Mollin, R.A., \textit{Polynomial solutions for Pell's equation revisited}, Indian Journal of Pure and Applied Mathematics \textbf{28} (1997), 429-438.
\bibitem{Nathanson} Nathanson, Melvyn B., \textit{Polynomial Pell's equations}, Proceedings of the American Mathematical Society \textbf{56} (1976), No. 1, 89-92.
\bibitem{Ramasamy} Ramasamy, Ams, \textit{Polynomial solutions for the Pell's equation}, Indian Journal of Pure and Applied Mathematics \textbf{25} (1994), 577-577.
\bibitem{Sage} SageMath, the Sage Mathematics Software System (Version 10.0),
   The Sage Developers, 2023, http://www.sagemath.org.
\bibitem{SageCode} Scherr, Z., Katherine, T. Computations For Square-free, Quartic, Monic Pellian Polynomials Over Z[x] [Computer software]. \url{https://github.com/zlscherr/QuarticIntegralMonicPell}
\bibitem{schmidt2000continued} Schmidt, Wolfgang, \textit{On Continued Fractions and Diophantine Approximation in Power Series Fields}, Acta Arithmetica \textbf{95} (2000), No. 2, 139-166.
\bibitem{WY1} Webb, William, Yokota, Hisashi, \textit{Polynomial Pell's equation}, Proceedings of the American Mathematical Society \textbf{131} (2003), No. 4, 993-1006.
\bibitem{WY2} Webb, William, Yokota, Hisashi, \textit{Polynomial Pell's equation--ii}, Journal of Number Theory \textbf{106} (2004), No. 1, 128-141.
\bibitem{Yokota1} Yokota, Hisashi, \textit{Polynomial Pell's equation and periods of quadratic irrationals}, JP Journal of Algebra, Number Theory, and Applications \textbf{8} (2007), 135-144.
\bibitem{Yokota2} Yokota, Hisashi, \textit{Solutions of polynomial Pell's equations}, Journal of Number Theory \textbf{130} (2010), No. 9, 2003-2010.
\end{thebibliography}
\end{document}